\documentclass[12pt]{amsart}
\usepackage{amsmath,amsthm,amscd,euscript,amsfonts}
\def \r{\mathbb R}
\def \q{\mathbb Q}
\def \z{\mathbb Z}

\theoremstyle{remark}

\theoremstyle{plane}

\numberwithin{equation}{section}

\title
[Two-dimensional continued fractions.]
{On two-dimensional continued fractions for the integer hyperbolic
matrices with small norm.}
\author{Oleg Karpenkov}
\thanks{This work is partially supported by INTAS-00-0259, NS-1972.2003.1 and
RFFI-01-01-00660.}
\date{UDK: 511.9}
\email[Oleg Karpenkov]{karpenk@mccme.ru}

\begin{document}
\input epsf
\maketitle

We say that two continued $(k-1)$-dimensional fractions (see~[1])
are {\it equivalent } if there exists integer lattice preserving
linear transform of the coordinates for the space $\r^{k}$
that maps one continued fraction to the other.
In this note we give the classification of two-dimensional continued
fractions constructed by the matrices whose norms are less than seven
up to the mentioned equivalence relation for multidimensional continued
fractions. This classification is based on the properties
for the family of Frobenius type matrices.
The norm here is the sum of the absolute values over all coefficients
of the matrix.
In conclusion an example of non-Frobenius matrix is constructed.
The author is grateful to V.~I.~Arnold for constant attention to this
work and useful remarks.

By $M(k,\z)$ we denote the set of $k\times k$ matrices with integer
coefficients and irreducible over the field $\q$ characteristic
polynomial. Note, that there is no matrices in the set $M(3,\z)$
with norm less than or equal to three. 
The numbers of matrices for the set $M(3,\z)$ with the norm
four, five, and six are equals to 240, 1248 and 8112 correspondingly.

Let $H(k,\z)$ be the subset of such matrices from $M(k,\z)$,
that it is possible to construct multidimensional continued fractions,
This means that $H(k,\z)$ is a subset of {\it hyperbolic} matrices
with all real and distinct eigenvalues.
The lowest possible norm for the matrices if the set $H(3,\z)$
is equal to five.
The number of matrices of the set $H(3,\z)$ with norm
five and six are equal to 48 and 912 correspondingly.

The $k \times k$-matrix obtained from the identity matrix of the order
$(k-1)\times (k-1)$ by adding the column of zeroes from the left
and the row of the numbers $(a_k,\ldots,a_1)$ from the bottom
is called {\it Frobenius} matrix and denoted by $M_{a_1,\ldots , a_k}$.
Note that the characteristic polynomial for such matrix is equal to
$\chi(x)=(-1)^k(x^k-a_1x^{k-1}-\cdots -a_{k-1}x-a_k)$. 
Furthermore the matrix $M_{a_1, \cdots, a_k}$ coincides with the
matrix of the multiplication by the element $x$ in the field
$\q[x]\big/ (\chi(x))$.

A matrix from the set $H(k,\z)$ is said to be of {\it Frobenius type}
if the multidimensional continued fraction for this matrix
is equivalent to the multidimensional continued fraction
constructed by some Frobenius matrix with integer coefficients.

Here we give the equivalent definition:
the matrix $A$ from $H(k,\z)$ is called the{\it Frobenius type} matrix
if there exist integer numbers $a_1,\ldots ,a_k$ and a matrix $X$
in the group $SL(k,\z)$ such that the matrix $XAX^{-1}$
commutes with the matrix $M_{a_1,\ldots ,a_k}$.
This definition can be naturally extended to the case of matrices 
for the set $M(k,\z)$:
the set $H(k,\z)$ should be replaced with the set $M(k,\z)$.

{\bf Theorem~1.}
All matrices with the norm less than or equal to six in the set $M(3,\z)$
are of Frobenius type. 

Theorem~1 essentially facilitates the classification of two-dimensional 
continued fractions constructed by matrices of cubic irrationalities with
norm no greater than six up to the multidimensional continued fraction
equivalence relation.

So, there is no matrices with the norms less than five
among the matrices of the set $H(3,\z)$. 
If the norm of the matrix from this set is equal to five,
then the two-dimensional continued fraction constructed by
this matrix is equivalent to the continued fraction constructed 
by the matrix $M_{-1,2,1}$.
This fraction is called {\it the golden ratio}.
The homeomorphic type of the torus decomposition
for the sail containing the point $(0,0,1)$ of this two-dimensional
continued fraction is the following: 
\begin{tabular}{c}
$\epsfbox{fig.1}.$ 
\end{tabular}
If the norm of the matrix is equal to six than the following
three cases are possible:
the corresponding continued fraction is equivalent either
to the golden ratio (for 480 matrices);
or to the continued fraction for the matrix $M_{-1,3,1}$:
\begin{tabular}{c}
$\epsfbox{fig.2}$ 
\end{tabular}
(for 240 matrices);
or to the continued fraction for the matrix $M_{0,3,1}$:
\begin{tabular}{c}
$\epsfbox{fig.3}$ 
\end{tabular}
(for 192 matrices).

The given classification and the proof of Theorem~1
are based on the generalization of the following theorem.

{\bf Theorem~2.}
The matrix $A=(a_{ij})$ from the set $M(2,\z)$, where $1\le i,j \le 2$,
is of Frobenius type iff
there exists the integer solution for the following equation:
$|Q_A(x,y)|=1$, where 
$Q_A(x,y)=\big( a_{12}x^2+(a_{22}-a_{11})xy-a_{21}y^2\big)
\big(\gcd(a_{12},a_{21},a_{22}-a_{11})\big)^{-1}$.

Let us make necessary notations before formulating the generalization 
of Theorem~2 for the case of matrices in the set $M(3,\z)$.

Denote by $\bar{\Xi} (С)$ the set of all integer matrices commuting with $C$.

{\bf Statement.}
The set of matrices $\bar{\Xi} (С)$ (for some matrix $С$ in the set $M(k,\z)$)
is the additive group that is isomorphic to $\z ^k$.
 
Let the unit matrix $E$ and integer matrices
$A$ and $B$ generates the basis of the group $\bar{\Xi}(C)$, where
$A=(a_{ij})$, $1\le i,j\le 3$. 
Since $B$ is contained in the group $\bar{\Xi}(C)$,
the matrix $B$ is equal to some polynomial of the matrix $A$,
and the power of this polynomial is no greater than two: 
$B=\alpha A^2 +\beta A+ \gamma E$,
here the coefficients $\alpha$, $\beta$ and $\gamma$ of the polynomial
are rational. 
Let the characteristic polynomial of the matrix $A$ equals
$\chi(x)=-x^3+a_1x^2-a_2x+a_3$. Let by definition
$$
\begin{array}{c}
\bar P(m,n)=m^3+\big( 2a_1\alpha +3\beta \big) m^2n+
\big( (a_2+a_1^2)\alpha^2+4a_1\alpha \beta+ 3\beta ^2 \big) mn^2+\\
\quad + \big( (a_1 a_2 -a_3)\alpha^3+(a_2 + a_1^2)\alpha^2\beta+
2a_1 \alpha \beta^2 + \beta ^3 \big)n^3.\\
\end{array}
$$

For two matrices $A=(a_{ij})$ and $B=(b_{ij})$ denote
$
\langle ij,kl\rangle _{A,B}=a_{ij}b_{kl}-a_{kl}b_{ij}$
(all indices are integer numbers from one to three).
Let us define the following polynomial:
$$
\begin{array}{l}
\tilde P_{A,B}(x,y,z)=
\langle 12,13\rangle _{A,B}x^3+
\langle 23,21\rangle _{A,B}y^3+
\langle 31,32\rangle _{A,B}z^3+\\
\big( \langle 13,11\rangle _{A,B}+
\langle 22,13\rangle _{A,B}+
\langle 12,23\rangle _{A,B} \big)x^2y+
\big( \langle 22,23\rangle _{A,B}+
\langle 23,11\rangle _{A,B}+
\langle 13,21\rangle _{A,B} \big)xy^2+\\
\big( \langle 33,12\rangle _{A,B}+
\langle 12,11\rangle _{A,B}+
\langle 32,13\rangle _{A,B} \big)x^2z+
\big( \langle 32,33\rangle _{A,B}+
\langle 11,32\rangle _{A,B}+
\langle 31,12\rangle _{A,B} \big)xz^2+\\
\big( \langle 21,22\rangle _{A,B}+
\langle 33,21\rangle _{A,B}+
\langle 23,31\rangle _{A,B} \big)y^2z+
\big( \langle 31,22\rangle _{A,B}+
\langle 33,31\rangle _{A,B}+
\langle 21,32\rangle _{A,B} \big)yz^2+\\
\big( \langle 11,22\rangle _{A,B}+
\langle 22,33\rangle _{A,B}+
\langle 33,11\rangle _{A,B}+
3\langle 13,31\rangle _{A,B} \big)xyz.
\end{array}
$$

{\bf Theorem~3.}
A matrix $C$ in the set $M(3,\z)$ is a Frobenius type matrix iff
there exist the integer solution of the following equation with
integer coefficients in the variables $x$, $y$, $z$, $m$ and $n$:
$|Q_{A,B}(x,y,z;m,n)|=1$, where 
$Q_{A,B}(x,y,z;m,n)=\bar P_{A,B}(m,n) \tilde P_{A,A^{\vee}}(x,y,z)$.
Here $E$, $A$ and $B$ generate the basis of the group $\bar{\Xi}(C)$,
and the matrix $A^{\vee}$ denotes such matrix that in the $ij$ place
of it the corresponding complementary minor multiplied by $(-1)^{i+j}$
is standing.

In conclusion we give the example of the hyperbolic matrix in the set
$M(3,\z) \cap SL(3,\z)$ that is not of Frobenius type.

{\bf The corollary of Theorem~3}.
The following matrix is not of the Frobenius type:

$$A=\left(
\begin{array}{ccc}
1 & 2 & 0\\
0 & 1 & 2\\
-7 & 0 & 29\\
\end{array}
\right).
$$

{\it Proof}.
Consider the polynomial $Q_{A,B}$ for the basis $E$, $A$
and $B=\frac{A^2}{2}-15A+\frac{29E}{2}$ of the class $\bar{\Xi}(A)$:
$
\big(
4x_1^3+56x_1^2x_2-14x_2^3+784x_1^2x_3+392x_1x_3^2-196x_2^2x_3+49x_3^3+42x_1x_2x_3
\big)\big(
2m^3-28m^2n+7n^3
\big).
$
The equation $|Q_{A,B}|=1$ has no integer solutions since
$2m^3-28m^2n+7n^3=\pm 1$ is not solvable modulo seven.
Therefore by Theorem~3 the matrix $A$  is not of Frobenius type.
\qed

The norm of the matrix cited above is $42$.
There is unknown if there exists the non-Frobenius matrix with a lesser norm.

\begin{center}
{\bf REFERENCES.}
\end{center}

1. Arnold~V.~I. Continued fractions. M.: MCCME, 2002.

\end{document}